\documentclass[a4paper,12pt]{amsart}

\usepackage{amssymb}
\usepackage{latexsym} 
\usepackage{amsfonts} 
\usepackage{amsmath}
\usepackage{eucal} 
\usepackage{bm} 
\usepackage{bbm} 
\usepackage{graphicx} 
\usepackage[english]{varioref} 
\usepackage[nice]{nicefrac} 
\usepackage[all]{xy}
\usepackage{amsthm}

\newcommand{\supp}{\text {\rm supp}}

\def\ge{\geqslant}
\def\le{\leqslant}
\def\a{\alpha}

\def\g{\gamma}

\def\d{\delta}
\def\D{\Delta}

\def\e{\epsilon}

\def\s{\sigma}

\def\th{\theta}

\def\l{\lambda}

\def\i{^{-1}}
\def\ZZ{\mathbb Z}
\def\NN{\mathbb N}
\def\QQ{\mathbb Q}
\def\FF{\mathbb F}

\newcommand{\kk}{\Bbbk}

\def\co{\mathcal O}

\def\tS{\tilde S}

\def\tW{\tilde W}
\def\tw{\tilde w}

\def\fo{\mathfrak o}

\def\<{\langle} 
\def\>{\rangle}

\theoremstyle{plain}
\newtheorem{thm}{Theorem}[section] 
\newtheorem*{thm*}{Theorem} 
 \newtheorem{prop}[thm]{Proposition}
 \newtheorem{lem}[thm]{Lemma}
 \newtheorem{cor}[thm]{Corollary}

\theoremstyle{definition}

\theoremstyle{remark}

\newtheorem*{rmk*}{Remark}
\newtheorem*{claim*}{Claim}

\begin{document}

\author{Xuhua He}
\address{Department of Mathematics, The Hong Kong University of Science and Technology, Clear Water Bay, Kowloon, Hong Kong}
\thanks{The author is partially supported by HKRGC grants 601409.}
\email{maxhhe@ust.hk}
\title[]{Closure of Steinberg fibers and affine Deligne-Lusztig varieties}

\begin{abstract}
We discuss some connections between the closure $\bar F$ of a Steinberg fiber in the wonderful compactification of an adjoint group and the affine Deligne-Lusztig varieties $X_w(1)$ in the affine flag variety. Among other things, we describe the emptiness/nonemptiness pattern of $X_w(1)$ if the translation part of $w$ is quasi-regular. As a by-product, we give a new proof of the explicit description of $\bar F$, first obtained in \cite{H1}. 
\end{abstract}

\maketitle

\section*{Introduction} 
\subsection{} Let $\kk$ be an algebraic closure of a finite field $\FF_q$, $L=\kk((\e))$ be the formal Laurent series and $\fo=\kk[[\e]]$ be the ring of formal power series. Let $\s$ be an automorphism on $L$ defined by $\s(\sum a_n \e^n)=\sum a_n^q \e^n$. 

Let $G$ be a simple algebraic group over $\kk$, split over $\FF_q$. We fix opposite Borel subgroups $B$ and $B^-$. Let $K=G(\fo)$ be a maximal bounded subgroup of the loop group $G(L)$. Let $I$ (resp. $I'$) be the inverse image of $B^-$ (resp. $B$) under the projection map $K \mapsto G$ sending $\e$ to $0$. The automorphism $\s$ on $L$ induces an automorphism on $G(L)$, which we still denote by $\s$. 

For any element $w$ in the extended affine Weyl group $\tW$ of $G(L)$, the affine Deligne-Lusztig variety $X_w(1)$ is defined by $$X_w(1)=\{g \in G(L)/I; g \i \s(g) \in I \dot w I\}.$$ 

\subsection{} One interesting question is to determine when $X_w(1)$ is empty. This is related to the reduction of Shimura varieties. See the survey papers by Rapoport \cite{Ra} and Haines \cite{Ha}. 

Reuman \cite{R} gave a conjecture on the emptiness/nonemptiness pattern of $X_w(1)$ if $w$ lies in the lowest two-sided cell of $\tW$ (i.e. the union of shrunken Weyl chambers). The emptiness is proved by G\"ortz, Haines, Kottwitz and Reuman in \cite[Proposition 9.5.4]{GHKR2} and the nonemptiness is proved in a joint work with G\"ortz in \cite{GH}. A partial result on the nonemptiness is also obtained by Beazley in \cite{Be}. 
 
For the elements outside the lowest two-sided cell, a conjecture is given in \cite[Conjecture 1.1.1]{GHKR2} and it is proved in \cite[Theorem 1.1.2]{GHKR2} that the emptiness occurs as predicted. But it is unknown if $X_w(1)$ is nonempty for the remaining elements $w$. 

\subsection{} In this paper, we will study the emptiness/nonemptiness pattern of $X_w(1)$ from a different point of view. We relate this problem to the problem of describing the closure $\bar F$ of Steinberg fiber in the wonderful compactification $\bar G$. The latter problem was solved in an earlier paper \cite{H1}. 

There is a specialization map from the loop group $G(L)$ to the wonderful compactification $\bar G$, introduced by Springer in \cite{Sp}. We will show in the paper that this map gives nice correspondences between the decomposition of $G(L)$ into $I' \times I'$-orbits and the decomposition of $\bar G$ into $B \times B$-orbits, and between the decomposition of $G(L)$ into $K$-stable pieces and the decomposition of $\bar G$ into $G$-stable pieces. Moreover, this map connects the emptiness/nonemptiness pattern of $X_w(1)$ to the emptiness/nonemptiness of the intersection of the closure of $\bar F$ with certain $B \times B$-orbit in $\bar G$. 

The main purpose of this paper is to prove that if the translation part of $w$ is quasi-regular (see $\S$\ref{qr} for the definition), then $X_w(1)$ is empty (resp. nonempty) if and only if the corresponding intersection in $\bar G$ is empty (resp. nonempty). This  includes a large fraction of the cases outside the lowest two-sided cell and the generic cases in the lowest two-sided cell. 

The precise statements are Proposition \ref{main1} (for the emptiness) and Theorem \ref{main2} (for the nonemptiness). Some special cases are stated in Corollary \ref{main3} without using the wonderful compactification. 

It would be interesting to connect the emptiness/nonemptiness pattern here with \cite[Conjecture 1.1.1]{GHKR2}. 

\subsection{} We now review the content of this paper in more detail.

In section 1, we recall the definition of the wonderful compactification and Springer's specialization map $G(L) \to \bar G$. 
The correspondence between the decomposition of $G(L)$ into $I' \times I'$-orbits and the decomposition of $\bar G$ into $B \times B$-orbits follows easily from the definition. In section 2, we discuss the correspondence between the decomposition of $G(L)$ into $K$-stable pieces and the decomposition of $\bar G$ into $G$-stable pieces. The closure relation of the $K$-stable pieces is also obtained. In section 3, we discuss some connections between the affine Deligne-Lusztig varieties and the closure of Steinberg fibers in $\bar G$. We also give a new proof of \cite[Theorem 4.3]{H1}. In section 4, we prove our main result on the emptiness/nonemptiness pattern for affine Weyl group element with quasi-regular translation part. Further discussions on the intersection of the closure of $\bar F$ with certain $B \times B$-orbit in $\bar G$ will be discussed in section 5. 

\section{Some decompositions on $\bar G$ and $G(L)$}

\subsection{} Let $B$ be a Borel subgroup of $G$, $B^-$ be an opposite Borel subgroup and $T=B \cap B^-$. Let $X$ be the coroot lattice and $Y$ be the coweight lattice. We denote by $Y^+$ the set of dominant coweights and $X^+=X \cap Y^+$. Let $(\a_i)_{i \in S}$ be the set of simple roots determined by $(B, T)$. Let $\Phi$ (resp. $\Phi^+$, $\Phi^-$) be the set of roots (resp. positive roots, negative roots). We denote by $W$ the Weyl group $N(T)/T$. For $i \in S$, we denote by $s_i$ the simple reflection corresponding to $i$. For $w \in W$, we denote by $\supp(w)$ the set of simple reflections occurring in a reduced expression of $w$. For $w \in W$, we choose a representative $\dot w$ in $N(T)$. 

For any $J \subset S$, let $P_J \supset B$ be the standard parabolic subgroup corresponding to $J$ and $P_J^- \supset B^-$ be the opposite parabolic subgroup. Let $L_J=P_J \cap P_J^-$ be a Levi subgroup. Let $\Phi_J$ be the roots of $L_J$, i.e., the roots spanned by $\a_j$ for $j \in J$. 

For any parabolic subgroup $P$, let $U_P$ be the unipotent radical of $P$. We simply write $U$ for $U_B$ and $U^-$ for $U_{B^-}$. 

For any $J \subset S$, let $\rho^\vee_J$ be the dominant coweight with $$\<\rho^\vee_J, \a_i\>=\begin{cases} 1, & \text{ if } i \in J \\ 0, & \text{ if } i \notin J\end{cases}.$$ We simply write $\rho^\vee$ for $\rho^\vee_S$. 

\subsection{} Let $\tW=N(T(L))/(T(L) \cap I)$ be the extended affine Weyl group of $G(L)$. It is known that $\tW=W \ltimes Y=\{w \e^\chi; w \in W, \chi \in Y\}$. We call $\chi$ the {\it translation part} of $w \e^\chi$. Let $l: \tW \to \NN \cup \{0\}$ be the length function. For $x=w \e^\chi \in \tW$, let $\dot x=\dot w \e^\chi$ be a representative in $N(T(L))$. 

Let $W_a=W \ltimes X \subset \tW$ be the affine Weyl group. Set $\tS=S \cup \{0\}$ and $s_0=\e^{\th^\vee} s_\th$, where $\th$ is the largest positive root of $G$. Then $(W_a, \tS)$ is a Coxeter system. For any $J \subset \tS$, let $W_J$ be the subgroup of $W_a$ generated by $J$ and $\tW^J$ be the set of minimal length coset representative of $\tW/W_J$. In the case where $J \subset S$, we write $W^J$ for $\tW^J \cap W$. 

For $w \in \tW^J$, set $$I(J, w)=\max\{K \subset J; \forall i \in K, \exists j \in K \text{ such that } w s_i=s_j w\}.$$

In particular, an element $w \in \tW^S$ is of the form $x \e^{-\l}$, where $\l \in Y^+$ and $x \in W^{I(\l)}$ and $I(S, x \e^{-\l})=I(I(\l), x)$, here $I(\l)=\{i \in S; \<\l, \a_i\>=0\}$. 

\subsection{}\label{bb} Let $\bar G$ be the wonderful compactification of $G$ (\cite{DP}, \cite{Str}). It is an irreducible, smooth projective $(G \times G)$-variety with finitely many $G \times G$-orbits $Z_J$ indexed by the subsets $J$ of $S$. We follow \cite[section 2]{Sp1} for the description of $Z_J$. 

Let $\l$ be a dominant coweight. Since $\bar G$ is complete, $\l(0)=\lim_{\e \to 0} \l(\e)$ is a well-defined point of $\bar G$. Moreover, if $\mu$ is another dominant coweight with $I(\l)=I(\mu)$, then $\l(0)=\mu(0)$. Set $h_J=\l(0)$ for dominant coweight $\l$ with $I(\l)=J$. This is a base point of $Z_J$. The map $(x, y) \mapsto (x, y) \cdot h_J$ induces an isomorphism from the quotient space $(G \times G) \times_{P^-_J \times P_J} L_J/Z(L_J)$ to $Z_J$. 

\subsection{}\label{gs}  Now we recall two partitions of $\bar G$. 

For $J \subset S$, $x \in W^J$ and $y \in W$, set $$[J, x, y]=(B \dot x, B \dot y) \cdot h_J.$$ Then $[J, x, y]$ is a $B \times B$-orbit in $\bar{G}$. By \cite{B1} and \cite[Lemma 1.3]{Sp}, \[\tag{1} \bar G=\sqcup_{J \subset S, x \in W^J, y \in W} [J, x, y].\]

Let $G_\D$ be the diagonal image of $G$ in $G \times G$. For $J \subset S$ and $w \in W^J$, set $$Z_{J, w}=G_\D \cdot [J, w, 1].$$ We call $Z_{J, w}$ a $G$-stable piece of $\bar G$. By \cite[12.3]{L1} and \cite[Prop 2.6]{H1}, 
\[\tag{2} \bar G=\sqcup_{J \subset S} \sqcup_{w \in W^J} Z_{J, w}.\]

The following properties will also be used in this paper. 

(3) For $J \subset S$ and $w \in W^J$, $\overline{Z_{J, w}}=\sqcup_{K \subset J, w' \in W^J, w' \ge u w u \i \text{ for some } u \in W_J} Z_{K, w'}$. See \cite[Theorem 4.5]{H3}. 

(4) Let $J \subset S$ and $w \in W^J$. If $w$ is a Coxeter element, then $G$ acts transitively on $Z_{J, w}$. See \cite[Corollary 3.8]{Sp3}.

\subsection{}\label{sp1} Since $\bar G$ is complete, the inclusion $\fo \to L$ induces a bijection from $\bar G(\fo)$ to $\bar G(L)$. Now the specialization $\e \mapsto 0$ defines a map $s: \bar G(L) \to \bar G$. In particular, for any $g \in G(L)$, $s(g) \in \bar G$. This is the specialization map introduced in \cite[2.1]{Sp2}. In particular, we have that $s(K)=G$ and $s(I')=B$.

Notice that any element in $\tW$ can be written in a unique way as $x \e^\l y \i$, where $\l \in Y^+$, $x \in W^{I(\l)}$ and $y \in W$. We have the following decompositions on $G(L)$, \[\tag{1} G(L)=\sqcup_{\l \in Y^+} K \e^\l K=\sqcup_{\l \in Y^+, x \in W^{I(\l)}, y \in W} I' \dot x \e^\l \dot y \i I'.\]

We have that \begin{gather*} s(K \e^\l K)=(s(K), s(K)) \cdot s(\e^\l)=(G, G) \cdot h_{I(\l)}=Z_{I(\l)}, \\ s(I' \dot x \e^\l \dot y \i I')=(s(I') \dot x, s(I') \dot y) \cdot h_{I(\l)}=[I(\l), x, y]. \end{gather*}

Thus by the decomposition (1) and the decomposition \ref{gs} (1), we have that 

(2) For any $J \subset S$, $s \i(Z_J)=\sqcup_{\l \in Y^+, I(\l)=J} K \e^\l K$.

(3) For any $J \subset S, x \in W^J, y \in W$, $s \i([J, x, y])=\sqcup_{\l \in Y^+, I(\l)=J} I' \dot x \e^\l \dot y \i I'$. 

\section{$K$-stable pieces in $G(L)$ and $G$-stable pieces in $\bar G$}

\subsection{} An analogue of $G$-stable pieces in the affine case is introduced in \cite{L2}. For any $w \in \tW^S$, set $$K_w=K \cdot I \dot w I, \qquad K_{w, \s}=K \cdot_\s I \dot w I,$$ here $\cdot$ is the usual conjugation action of $G(L)$, $g \cdot g'=g g' g \i$ and $\cdot_\s$ is the $\s$-twisted conjugation action of $G(L)$, $g \cdot_\s g'=g g' \s(g) \i$. We call $K_w$ a $K$-stable piece of $G(L)$ and $K_{w, \s}$ a $\s$-twisted $K$-stable piece of $G(L)$. 

This definition is different from the one in \cite{L2}. However, one can show in the same way as in \cite[Prop 2.6]{H1} that the two definitions are equivalent. Since the equivalence of the two definition is not used in this paper, we skip the details. 

\subsection{}\label{min} We first recall some properties about ``partial conjugation action'' of $W$ on $\tW$. Although the results were proved for affine Weyl groups in \cite{H4}, it is easy to see that they also hold for extended affine Weyl groups. 

Let $J \subset \tS$. We consider the conjugation action of $W$ on $\tW$, $x \cdot y=x y x \i$ for $x \in W$ and $y \in \tW$. For $w \in \tW^S$, set $$[w]=W \cdot (w W_{I(S, w)})=W \cdot (W_{I(S, w)} w).$$ By \cite[Corollary 2.6]{H4}, $\tW=\sqcup_{w \in \tW^S} [w]$. 

Given $w, w' \in \tW$ and $i \in S$, we write $w \xrightarrow{s_i} w'$ if $w'=s_i w s_i$ and $l(w') \le l(w)$. If $w=w_0, w_1, \cdots, w_n=w'$ is a sequence of elements in $\tW$ such that for all $k$, we have $w_{k-1} \xrightarrow{s_j} w_k$ for some $j \in S$, then we write $w \rightarrow_S w'$. We write $w \approx_S w'$ if $w \to_S w'$ and $w' \to_S w$.

By \cite[Proposition 3.4]{H4}, we have the following properties:

(1) for any $w \in \tW$, there exists a minimal length element $w' \in W \cdot w$ such that $w \to_S w'$. Moreover, we may take $w'$ to be an element of the form $v w_1$ with $w_1 \in \tW^S$ and $v \in W_{I(S, w_1)}$.

(2) Let $\co$ be a $W$-orbit of $\tW$ that contains $w \in \tW^S$. Then $w' \approx_S w$ for any minimal length element $w' \in \co$. 

\subsection{}\label{min1} By \cite[Corollary 4.5]{H4}, for any $W$-orbit $\co$ of $\tW$ and $v \in \co$, the following conditions are equivalent:

(1) $v$ is a minimal element in $\co$ with respect to the restriction to $\co$ of the Bruhat order on $\tW$.

(2) $v$ is an element of minimal length in $\co$.

We denote by $\co_{\min}$ the set of elements in $\co$ satisfy the above conditions. 

As in \cite[4.7]{H4}, we have a natural partial order $\le_S$ on $\tW^S$ defined as follows:

Let $w, w' \in \tW^S$. Then $w \le_S w'$ if for some (or equivalently, any) $v' \in (W \cdot w')_{\min}$, there exists $v \in (W \cdot w)_{\min}$ such that $v \le v'$.

In general, for $w \in \tW^S$ and $w' \in \tW$, we write $w \le_S w'$ if there exists $v \in (W \cdot w)_{\min}$ such that $v \le w'$. By \cite[Lemma 4.4]{H4}, if $w'' \to_S w'$ and $w \le_S w'$, then $w \le_S w''$. 

\

The main purpose of this section is to prove the following correspondence between $K$-stable pieces in $G(L)$ and $G$-stable pieces in $\bar G$. 

\begin{thm}\label{corr}
Let $*$ be the involution on $W$ defined by $w^*=w_0 w w_0$, where $w_0$ is the longest element in $W$. Then 

(1) For any $\l \in Y^+$ and $x \in W^{I(\l)}$, $s(K_{x \e^{-\l}})=Z_{I(-w_0 \l), x^*}$. 

(2) For any $J \subset S$ and $x \in W^J$, $s \i (Z_{J, x})=\sqcup_{\l \in Y^+, I(\l)=-w_0 J} K_{x^* \e^{-\l}}$.
\end{thm}

(1) We have that \begin{align*} s(K_{x \e^{-\l}}) &=s(K \cdot I \dot x \e^{-\l} I)=s(K \cdot \dot w_0 I \dot x \dot w_0 \e^{-w_0 \l} \dot w_0 I \dot w_0 \i) \\ &=s(K \cdot I' \dot x^* \e^{-w_0 \l} I')=s(K) \cdot s(I') \dot x^* s(\e^{-w_0 \l}) s(I') \\ &=G_\D (B \dot x^*, B) \cdot h_{I(-w_0 \l)}=Z_{I(-w_0 \l), x^*}.\end{align*}  

The proof of part (2) will be given in $\S$\ref{corr1}.

\

\begin{lem}\label{0} Let $w, w' \in \tW$ and $i \in S$ with $w'=s_i w s_i$. Then 

(1) If $l(w)=l(w')$, then $K \cdot I \dot w I=K \cdot I \dot w' I$.

(2) If $l(w')<l(w)$, then $K \cdot I \dot w I=K \cdot I \dot w' I \cup K \cdot I \dot s_i \dot w' I$
\end{lem}

If $w=w'$, then the statement is obvious. Now suppose that $w \neq w'$. By \cite[Lemma 1.6.4]{DL}, $l(w') \le l(w)$ implies that $s_i w<w$ or $w s_i<w$. 

If $s_i w<w$, then \[K \cdot I \dot w I=K \cdot I \dot s_i I \dot s_i \dot w I=K \cdot I \dot s_i \dot w I \dot s_i I=\begin{cases} K \cdot I \dot w' I, & \text{ if } l(w')=l(w) \\ K \cdot I \dot w' I \cup K \cdot I \dot s_i \dot w' I, & \text{ otherwise} \end{cases}.\] 

If $s_i w>w$ and $w s_i<w$, then $l(w')=l(w)$ and $K \cdot I \dot w I=K \cdot I \dot w \dot s_i I \dot s_i I=K \cdot I \dot s_i I \dot w \dot s_i I=K \cdot I \dot w' I$. \qed

\begin{lem}\label{1}
Let $w, w' \in \tW$. Then 

(1) If $w \to_S w'$, then $K \cdot I \dot w I \subset K \cdot I \dot w' I \cup \cup_{v \in w W_a, l(v)<l(w)} K \cdot I \dot v I$. 

(2) If $w \approx_S w'$, then $K \cdot I \dot w I=K \cdot I \dot w' I$. 
\end{lem}

By definition, there exists a finite sequence $w=w_0 \xrightarrow {i_1} w_1 \xrightarrow {i_2} \cdots \xrightarrow {i_m} w_n=w'$, where $i_j \in S$ for all $j$. We prove the lemma by induction on $m$. 

For $m=0$ this is clear. Now assume that $m>0$ and the statements hold for $m-1$. By the previous lemma, $K \cdot I \dot w I \subset K \cdot I \dot w_1 I \cup \cup_{x \in w W_a, l(x)<l(w)} K \cdot I \dot x I$. Notice that $l(w_1) \le l(w)$. By induction hypothesis, $K \cdot I \dot w_1 I \subset K \cdot I \dot w' I \cup \cup_{x \in w W_a, l(x)<l(w)} K \cdot I \dot x I$. Hence $K \cdot I \dot w I \subset K \cdot I \dot w' I \cup \cup_{x \in w W_a, l(x)<l(w)} K \cdot I \dot x I$. 

If moreover, $w \approx_S w'$, then $l(w_1)=l(w)$ and $w_1 \approx_S w'$. By induction hypothesis, $K \cdot I \dot w I=K \cdot I \dot w_1 I=K \cdot I \dot w' I$. \qed

\begin{lem}\label{2}
Let $w \in \tW^S$ and $v \in W_{I(S, w)}$. Then $K \cdot I \dot v \dot w I \subset K_w$. 
\end{lem}

Set $J=I(S, w)$. Notice that $I=(B^- \cap L_J) I_J$, where $I_J$ is the inverse image of $U_{P_J^-}$ under the map $K \to G$. Then $I_J$ is normal in $I$ and $\dot w (B^- \cap L_J) \dot w \i=B^- \cap L_J$. We have that \begin{align*} I \dot v \dot w I &=I_J (B^- \cap L_J) \dot v \dot w (B^- \cap L_J) I_J=I_J (B^- \cap L_J) \dot v (B^- \cap L_J) \dot w I_J \\ & \subset I_J L_J \dot w I_J.\end{align*} The map $l \mapsto \dot w l \dot w \i$ is an automorphism on $L_J$ and $\dot w (B^- \cap L_J) \dot w \i=B^- \cap L_J$. By \cite[Lemma 7.3]{Ste}, $L_J \dot w=\{l l' \dot w l \i; l \in L_J, l' \in B^- \cap L_J\}$. So \begin{align*} L_J \cdot I \dot w I_J &=\{l I_J (B^- \cap L_J) \dot w I_J l \i; l \in L_J\} \\ &=\{I_J l (B^- \cap L_J) \dot w l \i I_J; l \in L_J\}=I_J L_J \dot w I_J \end{align*} and $K \cdot I \dot v \dot w I \subset K \cdot L_J I_J \dot w I_J=K \cdot I \dot w I_J \subset K \cdot I \dot w I=K_w$. \qed

\begin{prop}\label{3}
Let $w \in \tW$. Then $$\overline{K \cdot I \dot w I}=K \cdot \overline{I \dot w I}=\cup_{w' \in \tW^S, w' \le_S w} K_w.$$ 
\end{prop}

Define the action of $I$ on $K \times G(L)$ by $i \cdot (k, g)=(k i \i, i g i \i)$. Let $K \times_I G(L)$ be the quotient space. Then the map $K \times G(L) \to G(L)$, $(k, g) \mapsto k g k \i$ induces a map $f: K \times_I G(L) \to G(L)$. Notice that each fiber is isomorphic to $K/I \cong G/B$. Thus $K \cdot \overline{I \dot w I}=f(K \times_I \overline{I \dot w I})$ is closed in $G(L)$. Hence $\overline{K \cdot I \dot w I}=K \cdot \overline{I \dot w I}$. 

If $w' \in W^S$ with $w' \le_S w$, then there exists an element $w'' \in (W \cdot w')_{\min}$ such that $w'' \le w$. By \ref{min} (2), $w'' \approx w'$. So by Lemma \ref{1}, $K_{w'}=K \cdot I \dot w'' I \subset K \cdot \overline{I \dot w I}$. Now we prove by induction on $l(w)$ that 

(a) $K \cdot I \dot w I \subset \cup_{w' \in \tW^S, w' \le_S w} K_{w'}$. 

If $w \in (W \cdot w)_{\min}$, then by \ref{min} (1), $w \approx v w_1$ for some $w_1 \in \tW^S$ and $v \in W_{I(S, w_1)}$. Then $w_1 \le_S w$. By Lemma \ref{2}, $K \cdot I \dot w I \subset K_{w_1}$. If $w \notin (W \cdot w)_{\min}$, then by \ref{min} (1), there exists $w_1 \approx w$ and $i \in S$ with $w_1 \xrightarrow i s_i w_1 s_i$ and $l(s_i w_1 s_i)<l(w)$. By Lemma \ref{0} and Lemma \ref{1}, $K \cdot I \dot w I=K \cdot I \dot w_1 I \subset K \cdot I \dot s_i \dot w I \cup K \cdot I \dot s_i \dot w \dot s_i I$. Since $l(s_i w_1 s_i)<l(w)$, we have that $s_i w_1 s_i<s_i w_1<w_1$. Thus for any $w' \in W^S$, if $w' \le_S s_i w_1 s_i$ or $w' \le_S s_i w_1$, then $w' \le_S w_1$ and $w' \le_S w$. (a) follows from induction hypothesis. 
 
Now $K \cdot \overline{I \dot w I}=\cup_{x \le w} K \cdot I \dot x I \subset \cup_{x \le w} \cup_{w' \in \tW^S, w' \le_S x} K_{w'}=\cup_{w' \in \tW^S, w' \le_S w} K_{w'}$. The proposition is proved. \qed

\begin{prop}\label{k}
We have that $G(L)=\sqcup_{w \in \tW^S} K_w$.
\end{prop}

\begin{rmk*}
This is essentially contained in \cite[1.4]{L2}. Here we give a different proof. 
\end{rmk*}

By Proposition \ref{3},  \begin{align*} G(L) &=\cup_{x \in \tW} I \dot x I=\cup_{x \in \tW} K \cdot I \dot x I=\cup_{x \in \tW} K \cdot \overline{I \dot x I} \\ &=\cup_{x \in \tW} \cup_{w \in \tW^S, w \le_S x} K_w=\cup_{w \in \tW^S} K_w. \end{align*} 

Let $w, w' \in \tW^S$. Then $w=x \e^{-\l}$ and $w'=x' \e^{-\l'}$ for some $\l, \l' \in Y^+$ and $x \in W^{I(\l)}, x' \in W^{I(\l')}$. Suppose that $K_w \cap K_{w'} \neq \emptyset$. Notice that $K_w \subset K \e^{-\l} K$ and $K_{w'} \subset K \e^{-\l'} K$. By the decomposition $G(L)=\sqcup_{\mu \in Y^+} K \e^{-\mu} K$, we have that $\l=\l'$. We also have that $s(K_w) \cap s(K_{w'}) \neq \emptyset$. Since $s(K_w)=Z_{-w_0 I(\l), x^*}$ and $s(K_{w'})=Z_{-w_0 I(\l), (x')^*}$, by the $G$-stable-piece decomposition \ref{gs} (2), we must have that $x=x'$.  Thus $w=w'$. The Proposition is proved. \qed

\begin{cor}
Let $w \in \tW^S$. Then $\overline{K_w}=\sqcup_{w' \in \tW^S, w' \le_S w} K_{w'}$. 
\end{cor}

\subsection{Proof of Theorem \ref{corr} (2)}\label{corr1}
By part (1), $K_{x^* \e^{-\l}} \subset s \i (Z_{J, x})$ for any $\l \in Y^+$ with $l(\l)=-w_0 J$. Now let $g \in s \i(Z_{J, x})$, then by Proposition \ref{k}, $g \in K_w$ for some $w \in \tW^S$. Again by part (1), $s(K_w)$ is a $G$-stable piece in $\bar{G}$. Hence $s(K_w)=Z_{J, x}$ and $w$ must be of the form $x^* \e^{-\l}$ for any $\l \in Y^+$ with $l(\l)=-w_0 J$. 

\subsection{} All the above results remain valid for the $\s$-twisted $K$-stable pieces of $G(L)$. In fact, we have stronger result below for $K_{w, \s}$ than Lemma \ref{2} for $K_w$. 

\begin{lem}\label{22}
Let $w \in \tW^S$ and $v \in W_{I(S, w)}$. Then $K \cdot_\s I \dot v \dot w I=K_{w, \s}$. 
\end{lem}

Set $J=I(S, w)$. By the proof of Lemma \ref{2}, $I \dot v \dot w I=I_J (B^- \cap L_J) \dot v (B^- \cap L_J) \dot w I_J$.  By Lang's theorem, $L_J=\{l \dot w \s(l) \i \dot w \i; l \in L_J\}$ and $L_J \dot w=\{l \dot w \s(l) \i\}$. Therefore $L_J \cdot_\s I \dot v \dot w I=L_J \cdot_\s I_J (B^- \cap L_J) \dot v (B^- \cap L_J) \dot w I_J=L_J \cdot_\s I_J \dot w I_J$. Taking $v=1$, then $L_J \cdot_\s I \dot w I=L_J \cdot_\s I_J \dot w I_J$. Hence $L_J \cdot_\s I \dot v \dot w I=L_J \cdot_\s I \dot w I$ and $K \cdot_\s I \dot v \dot w I=K_{w, \s}$.  \qed

\begin{lem}\label{33}
Let $w \in \tW$. Then $K \cdot_\s I \dot w I$ is a union of $\s$-twisted $K$-stable pieces. 
\end{lem}

We prove by induction on $l(w)$. 

If $w \in (W \cdot w)_{\min}$, then by \ref{min} (1), $w \approx  v w_1$ for some $w_1 \in \tW^S$ and $v \in W_{I(S, w_1)}$. Then $K \cdot_\s I \dot w I=K \cdot_\s I \dot v \dot w_1 I=K_{w_1, \s}$ is a $\s$-twisted $K$-stable piece.  If $w \notin (W \cdot w)_{\min}$, then by \ref{min} (1), there exists $w_1 \approx w$ and $i \in S$ with $w_1 \xrightarrow i s_i w_1 s_i$ and $l(s_i w_1 s_i)<l(w)$. Then $K \cdot_\s I \dot w I=K \cdot_\s I \dot w_1 I=K \cdot_\s I \dot s_i \dot w I \cup K \cdot_\s I \dot s_i \dot w \dot s_i I$. Since $l(s_i w_1 s_i), l(s_i w)<l(w)$, $K \cdot_\s I \dot s_i \dot w I$ and $K \cdot_\s I \dot s_i \dot w \dot s_i I$ are unions of $\s$-twisted $K$-stable pieces. Hence $K \cdot_\s I \dot w I$ is a union of $\s$-twisted $K$-stable pieces. \qed

\subsection{} By the same argument as in \cite[Corollary 2.6]{H5}, the subset $G_\s (B \dot x, B) \cdot h_J$ of $\bar G$ is a single $G_\s$-orbit, here $G_\s=\{(g, \s(g)); g \in G\}$, $J \subset S$ and $x \in W^J$. Thus using the specialization map $s: G(L) \to \bar G$, one can show that $K \cdot_\s I \dot w' I$ is a single orbit of $K_\s (U_K \times U_K)$ for any $w' \in \tW^S$, here $K_\s=\{(g, \s(g)); g \in K\} \subset K \times K$ and $U_K$ is the inverse image of $1 \in G$ under the projection map $K \to G(L)$ sending $\e$ to $0$. This gives another proof of the above Lemma. We omit the details. 

\section{ADLV and the closure of Steinberg fibers in $\bar{G}$}

We first discuss some equivalence conditions for the nonemptiness of an affine Deligne-Lusztig variety $X_w(1)=\{g \in G(L)/I; g \i \s(g) \in I \dot w I\}$. 

\begin{prop}\label{re}
Let $w \in \tW$ and $t \in I \cap T(L)$. Then the following conditions are equivalent:

(1) $X_w(1) \neq \emptyset$;

(2) $U(L) \cap K \cdot I \dot w I \neq \emptyset$;

(3) $t U(L) \cap K \cdot I \dot w I \neq \emptyset$;

(4) $\dot x \i U(L) \dot x \cap I \dot w I \neq \emptyset$ for some $x \in W$. 
\end{prop}

\begin{rmk*}
The equivalence of (1) and (4) is essentially contained in \cite[Section 6]{GHKR1}. 
\end{rmk*}

By definition, $X_w(1) \neq \emptyset$ if and only if $g \i \s(g) \in I \dot w I$ for some $g \in G(L)$. Using the decomposition $G(L)=\sqcup_{v \in \tW} U(L) \dot v I$, this is equivalent to $\dot v \i u \i \s(u) \s(\dot v) \in I \dot w I$ for some $u \in U(L)$ and $v \in \tW$. By \cite[1.3]{L0}, $\{u \i \s(u); u \in U(L)\}=U(L)$. Notice that $\s(\dot v) \in \dot v T$. Hence $1 \in G(L) \cdot_\s I \dot w I$ if and only if $\dot v \i U(L) \dot v \cap I \dot w I \neq \emptyset$ for some $v \in \tW$. Write $v$ as $v=\e^\l x$ for $x \in W$ and $\l \in Y$. Then $\dot v \i U(L) \dot v=\dot x \i \e^{-\l} U(L) \e^\l \dot x=\dot x \i U(L) \dot x$. Thus (1) is equivalent to (4). 

By the same argument, $t U(L) \cap G(L) \cdot I \dot w I \neq \emptyset$ if and only if $\dot x \i t U(L) \dot x \cap I \dot w I \neq \emptyset$ for some $x \in W$. Since $t \in I \cap T(L)$, this is equivalent to the condition (4). On the other hand, (4) implies (3) and (3) implies that $t U(L) \cap G(L) \cdot I \dot w I \neq \emptyset$. Hence (3) is equivalent to (4). Taking $t=1$, we obtain the equivalence between (2) and (4). \qed

\begin{cor}\label{inv}
Let $w \in \tW$. Then $X_w(1) \neq \emptyset$ if and only if $X_{w \i}(1) \neq \emptyset$.
\end{cor}

By Proposition \ref{re}, if $X_w(1) \neq \emptyset$, then $U(L) \cap K \cdot I \dot w I \neq \emptyset$. Applying the inverse map, we have that $U(L) \cap K \cdot I \dot w \i I \neq \emptyset$. Again by Proposition \ref{re}, $X_{w \i}(1) \neq \emptyset$. \qed 

\subsection{} Let $St: G \to T/W$ be the Steinberg map, i.e., $St(g)$ is the $W$-orbit in $T$ that contains an element conjugate to the semisimple part of $g$. The fibers of this map are called Steinberg fibers. It is known that each Steinberg fiber is of the form $G \cdot t U$ for some $t \in T$. An example of Steinberg fiber is the unipotent variety $\mathcal U$ of $G$. Some other examples are regular semisimple conjugacy classes of $G$. 

Let $F$ be a Steinberg fiber and $\bar F$ be its closure in $\bar G$. The following description of $\bar F$ was first obtained in \cite[Theorem 4.3 \& 4.5]{H1} using a case-by-case check. A more conceptual proof was obtained later in \cite{HT}. Springer also gave a different proof in \cite{Sp1}. Both \cite{HT} and \cite{Sp1} uses some properties of proper intersection. Below we give a group-theoretic proof in positive characteristic based on the connection between loop groups and group compactifications. 

\begin{thm}\label{uni}
Let $F$ be a Steinberg fiber of $G$. Then $$\bar F-F=\sqcup_{J \neq S} \sqcup_{w \in W^J, \supp(w)=S} Z_{J, w}.$$
\end{thm}

The part $\bar F-F \subset \sqcup_{J \neq S} \sqcup_{w \in W^J, \supp(w)=S} Z_{J, w}$ follows from the fact that the elements in $\bar F-F$ are represented by nilpotent endomorphisms. We refer to \cite[3.3(b)]{Sp1} for the details. 

By \ref{gs} (3), in order to show that $\sqcup_{J \neq S} \sqcup_{w \in W^J, \supp(w)=S} Z_{J, w} \subset \bar F$, it suffices to show that $Z_{J, w} \subset \bar F$ for any $J \neq S$ and Coxeter element $w \in W^J$. 

Assume that $F=G_\D \cdot t U$ for some $t \in T$. By \cite[1.4]{Sp2}, $\bar F=G_\D \cdot t \bar U$. Let $J \neq S$ and $w \in W^J$ be a Coxeter element of $W$. By \ref{gs} (4), $Z_{J, w}$ is a single $G$-orbit. So it suffices to show that $Z_{J, w} \cap t \bar U \neq \emptyset$. 

Pick $\l \in X^+$ with $I(-w_0 \l)=J$. Since $w$ and $w^*$ are Coxeter elements, there exists $\mu \in Y \otimes_\ZZ \QQ$ such that $\mu-(w^*) \i \mu=\l$. Let $m \in \NN$ with $m \mu \in Y$. Then $m \mu-(w^*) \i (m \mu)=m \l$ and $w^* \e^{-m \l}$ is conjugated to $w^*$ by $\e^{m \mu}$. Therefore $\dot w^* \in G(L) \cdot_\s \dot w^* \e^{-m \l}$. By Lang's theorem, $1 \in K \cdot_\s \dot w^*$. Hence $1 \in G(L) \cdot_\s I \dot w^* \e^{-m \l} I$. By Proposition \ref{re}, $t U(L) \cap K_{w^* \e^{-m \l}} \neq \emptyset$. Applying the specialization map $s: G(L) \to \bar G$, we have that $t \bar U \cap Z_{J, w} \neq \emptyset$. This finishes the proof. \qed

\

The following result relates the emptiness problem of ADLV to the description of the closure of a Steinberg fiber in $\bar G$. 

\begin{prop}\label{main1}
Let $J \subsetneqq S$, $x \in W^J$ and $y \in W$. If $[J, x, y] \subset \sqcup_{w \in W^J, \supp(w) \neq S} Z_{J, w}$, then for any $\l \in Y^+$ with $I(\l)=J$, we have that $X_{x \e^{-\l} y \i}(1)=\emptyset$.
\end{prop}

Let $\d$ be the diagram automorphism on $G$ whose induced permutation on $S$ equals $-w_0$. Then $\d$ also gives an automorphism on $G(L)$ which we still denote by $\d$. Then $\d(1)=1$ and $\d(I)=I$. So $\d(I \dot x \e^{-\l} \dot y \i I)=I \dot x^* \e^{w_0 \l} \dot (y^*) \i I$.  

If $X_{x \e^{-\l} y \i}(1) \neq \emptyset$, then $X_{x^* \e^{w_0 \l} (y^*) \i}(1) \neq \emptyset$ and $1 \in G(L) \cdot_\s I \dot x^* \e^{w_0 \l} (\dot y^*) \i I$. By Corollary \ref{re},  $U(L) \cap K \cdot I \dot x^* \e^{w_0 \l} (\dot y^*) \i I \neq \emptyset$. As in $\S$\ref{corr1}, $K \cdot I \dot x^* \e^{w_0 \l} (\dot y^*) \i I=K \cdot I' \dot x \e^{\l} \dot y \i I'$. Thus, $K \cdot U(L) \cap I' \dot x \e^{\l} \dot y \i I \neq \emptyset$. Applying the specialization map, we have that $G_\D \cdot \bar U \cap [J, x, y] \neq \emptyset$. That is impossible by the previous theorem.  \qed 

\section{Main result}

\subsection{}\label{qr} Note that $\sqcup_{w \in W_a} I \dot w I$ is a normal subgroup of $G(L)$ that contains $1$. It is easy to see that if $X_w(1) \neq \emptyset$ for some $w \in \tW$, then we must have that $w \in W_a$, i.e., the translation part of $w$ is $X$. 

Now let $\l \in X$. We call $\l$ {\it quasi-regular} if for any $\a \in \Phi$, either $\<\l, \a\>=0$ or $|\<\l, \a\>| \ge (\<\rho^\vee, \th\>+2)^{|S|+1}$. Any affine Weyl group element with quasi-regular translation part is of the form $x \e^{-\l} y \i$, where $\l \in X^+$ with $\<\l, \a_i\> \ge (\<\rho^\vee, \th\>+2)^{|S|+1}$ for any $i \notin I(\l)$, $x \in W^{I(\l)}$ and $y \in W$. Our main result below describes the emptiness/nonemptiness pattern of $X_w(1)$ if the translation part of $w$ is quasi-regular. 

\begin{thm}\label{main2}
Let $J \subsetneqq S$, $x \in W^J$, $y \in W$ and $\l \in X^+$ with $I(\l)=J$. Assume that $\l$ is quasi-regular. Then $X_{x \e^{-\l} y \i}(1) \neq \emptyset$ if and only if $[J, x, y] \nsubseteq \sqcup_{w \in W^J, \supp(w) \neq S} Z_{J, w}$.
\end{thm}

The ``only if'' is proved in Proposition \ref{main1}. The proof of ``if'' part will be given in $\S$\ref{main22}. Our strategy is as follows. First, we use the Proposition \ref{key1} below to reduce to problem to elements in $\tW^S$, using the technique of ``partial conjugation action'' introduced in \cite{H4}. Then in Proposition \ref{key2} we reduce the elements in $\tW^S$ with quasi-regular translation part to some elements for which the nonemptiness is already known. The trick we use here is similar to the ``P-operators'' introduced in \cite{H6}. Because of the quasi-regular condition, the case here is easier to handle than in loc. cit. 

\begin{prop}\label{key1}
Let $\l \in Y^+$ with $I(\l)=J \subset S$. Let $x, w \in W^J$ and $y \in W$. Then the following conditions are equivalent:

(1) $[J, x, y] \cap Z_{J, w} \neq \emptyset$;

(2) $I \dot x \e^{-\l} \dot y \i I \cap K_{w \e^{-\l}} \neq \emptyset$;

(3) $I \dot x \e^{-\l} \dot y \i I \cap K_{w \e^{-\l}, \s} \neq \emptyset$;

(4) $K_{w \e^{-\l}, \s} \subset K \cdot_\s I \dot x \e^{-\l} \dot y \i I$. 
\end{prop}

We first prove the equivalence of (1) and (2). 

Let $\d$ be the automorphism of $G(L)$ defined in the proof of Proposition \ref{main1}. Then \begin{gather*} G_\D \cdot [J, x, y]=s(K \cdot I' \dot x \e^\l \dot y \i I')=s(K \cdot I \dot x^* \e^{w_0 \l} (\dot y^*) \i I)=s(\d(K \cdot I \dot x \e^{-\l} \dot y \i I)), \\ Z_{J, w}=G_\D \cdot [J, w, 1]=s(\d(K \cdot I \dot w \e^{-\l} I))=s(\d(K_w)). \end{gather*} So if 
$I \dot x \e^{-\l} \dot y \i I \cap K_{w \e^{-\l}} \neq \emptyset$, then $K \cdot I \dot x \e^{-\l} \dot y \i I \cap K_{w \e^{-\l}} \neq \emptyset$ and $s(\d(K \cdot I \dot x \e^{-\l} \dot y \i I)) \cap s(\d(K_{w \e^{-\l}})) \neq \emptyset$. Thus $G_\D \cdot [J, x, y] \cap Z_{J, w} \neq \emptyset$ and $[J, x, y] \cap Z_{J, w} \neq \emptyset$. On the other hand, if $I \dot x \e^{-\l} \dot y \i I \cap K_{w \e^{-\l}}=\emptyset$, then by Proposition \ref{k}, $K \cdot I \dot x \e^{-\l} \dot y \i I \subset \sqcup_{w' \in W^J, w' \neq w} K_{w' \e^{-\l}}$. Hence \begin{align*} G_\D \cdot [J, x, y] &=s(\d(K \cdot I \dot x \e^{-\l} \dot y \i I)) \subset \sqcup_{w' \in W^J, w' \neq w} s(\d(K_{w' \e^{-\l}}) \\ &=\sqcup_{w' \neq W^J, w' \neq w} Z_{J, w}.\end{align*} Therefore $[J, x, y] \cap Z_{J, w}=\emptyset$. 

By Lemma \ref{33}, $K \cdot_\s I \dot x \e^{-\l} \dot y \i I$ is a union of $\s$-twisted $K$-stable pieces. Hence (3) is equivalent to (4). 

Set $\tw=x \e^{-\l} y \i$. Now we prove the equivalence of (2) and (3) by induction on $l(\tw)$. 

If $\tw \in (W \cdot \tw)_{\min}$, then by \ref{min} (1), $\tw \approx  v \tw_1$ for some $\tw_1 \in \tW^S$ and $v \in W_{I(S, \tw_1)}$. By Lemma \ref{2} and Lemma \ref{22}, $K \cdot I \dot \tw I \subset K_{\tw_1}$ and $K \cdot_\s I \dot \tw I=K_{\tw_1, \s}$. By Proposition \ref{k}, the $K$-stable pieces form a disjoint union of $G(L)$. Hence (2) is equivalent to $w \e^{-\l}=\tw_1$. Similarly,  (3) is equivalent to $w \e^{-\l}=\tw_1$. 

If $\tw \notin (W \cdot \tw)_{\min}$, then by \ref{min} (1), there exists $\tw_1 \approx \tw$ and $i \in S$ with $\tw_1 \xrightarrow i s_i \tw_1 s_i$ and $l(s_i \tw_1 s_i)<l(\tw)$. By Lemma \ref{0}, $K \cdot I \dot \tw I=K \cdot I \dot \tw_1 I=K \cdot I \dot s_i \dot \tw_1 I \cup K \cdot I \dot s_i \dot \tw_1 \dot s_i I$. Similarly, $K \cdot_\s I \dot \tw I=K \cdot_\s I \dot \tw_1 I=K \cdot_\s I \dot s_i \dot \tw_1 I \cup K \cdot_\s I \dot s_i \dot \tw_1 \dot s_i I$. Since $l(s_i \tw_1 s_i), l(s_i \tw)<l(\tw)$, by induction hypothesis, $I \dot s_i \dot \tw_1 I \cap K_{w \e^{-\l}} \neq \emptyset$ if and only if  $I \dot s_i \dot \tw_1 I \cap K_{w \e^{-\l}, \s} \neq \emptyset$ and $I \dot s_i \dot \tw_1 \dot s_i I \cap K_{w \e^{-\l}} \neq \emptyset$ if and only if  $I \dot s_i \dot \tw_1 \dot s_i I \cap K_{w \e^{-\l}, \s} \neq \emptyset$. So $K \cdot I \dot \tw I \cap K_{w \e^{-\l}} \neq \emptyset$ if and only if $K \cdot_\s I \dot \tw I \cap K_{w \e^{-\l}, \s} \neq \emptyset$. The equivalence of (2) and (3) is proved.  \qed

\begin{lem}
Let $x, y \in W$. Assume that for any $\a \in \Phi^+$ with $y \i \a \in \Phi^-$, $x \a \in \Phi^-$. Then $l(x y)=l(x)-l(y)$. 
\end{lem}

We prove by induction on $l(y)$. For $y=1$ this is clear. Suppose that $l(y) \ge 1$. Then there exists $i \in S$ such that $y \i \a_i \in \Phi^-$. Thus $y=s_i y'$ for some $y' \in W$ with $l(y')=l(y)-1$. We have that $x \a_i \in \Phi^-$. Hence $x=x' s_i$ for some $x' \in W$ with $l(x')=l(x)-1$. 

Now for any $\a \in \Phi^+$ with $y \i \a \in \Phi^-$, either $\a=\a_i$ or $s_i \a \in \Phi^+$ with $(y') \i (s_i \a) \in \Phi^-$. By our assumption, $x' (s_i \a)=x \a \in \Phi^-$. Therefore for any $\a' \in \Phi^+$ with $(y') \i \a' \in \Phi^-$, we have that $\a' \neq \a_i$ and $x' \a' \in \Phi^-$. By induction hypothesis, $l(x y)=l(x' y')=l(x')-l(y')=l(x)-l(y)$. \qed

\begin{prop}\label{key2}
Let $J \subset S$. Then for any $x \in W^J$ and $y \in W$ with $y \i x \in W^J$ and $\supp(y \i x)=S$ and any $\l \in X^+$ with $I(\l)=J$ and $\<\l, \a_i\> \ge (\<\rho^\vee, \th\>+2)^{|J|+1}$ for any $i \notin J$, we have that $X_{x \e^{-\l} y \i}(1) \neq \emptyset$. 
\end{prop}

We prove by induction on $|J|$. Note that \begin{align*} l(x \e^{-\l} y \i) &=l(\e^\l)-l(x)+l(y)=l(x \e^{-M \rho^\vee_{S-J}})+l(\e^{-\l+M \rho^\vee_{S-J}} y \i),\end{align*} where $M=(\<\rho^\vee, \th\>+2)^{|J|}$. So \begin{align*} \tag{a} K \cdot_\s I \dot x \e^{-\l} \dot y \i I &=K \cdot_\s I \dot x \e^{-M \rho^\vee_{S-J}} I \e^{-\l+M \rho^\vee_{S-J}} \dot y \i I\\ &=K \cdot_\s I \e^{-\l+M \rho^\vee_{S-J}} \dot y \i I \dot x \e^{-M \rho^\vee_{S-J}} I \\ & \supset K \cdot_\s I 
 \e^{-\l+M \rho^\vee_{S-J}} \dot y \i \dot x \e^{-M \rho^\vee_{S-J}} I \\ &=K \cdot_\s I \e^{-(\l-M \rho^\vee_{S-J}+y \i x M \rho^\vee_{S-J})} \dot y \i \dot x I. \end{align*}
 
Now $\l-M \rho^\vee_{S-J}+y \i x M \rho^\vee_{S-J}=v \g$ for some $\g \in X^+$ and $v \in W^{I(\g)}$. If $v=1$, then $\e^{-(\l-M \rho^\vee_{S-J}+y \i x M \rho^\vee_{S-J})} y \i x=\e^{-\g} y \i x$. By \cite{GH}, $X_{x \i y \e^\g}(1) \neq \emptyset$. Hence by Corollary \ref{inv}, $X_{\e^{-\g} y \i x}(1) \neq \emptyset$. By (a), $1 \in K \cdot_\s I \e^{-\g} \dot y \i \dot x I \subset K \cdot_\s I \dot x \e^{-\l} \dot y \i I$. Thus $X_{x \e^{-\l} y \i}(1) \neq \emptyset$ and the Proposition holds in this case. 

Now we consider the case where $v \neq 1$. For any $\a \in \Phi^+-\Phi^+_J$, \begin{align*} \tag{b} \<v \g, \a\> &=\<\l-M \rho^\vee_{S-J}, \a\>+\<y \i x M \rho^\vee_{S-J}, \a\>  \\ &  \ge  (\<\rho^\vee, \th\>+2)^{|J|+1}-M+M\<\rho^\vee_{S-J}, x \i y \a\> \\ & \ge  (\<\rho^\vee, \th\>+2)^{|J|+1}-M-M\<\rho^\vee, \th\>=M. \end{align*}

Let $\a \in \Phi^+$ with $v \i \a \in \Phi^-$. Then $\<v \g, \a\>=\<\g, v \i \a\> \le 0$. However, if $\<v \g, \a\>=0$, then $v \i \a \in \Phi^-_{I(\g)}$ and $\a=v (v \i \a) \in \Phi^-$. That is a contradiction. Hence $$\<v \g, \a\>=\<\l-M \rho^\vee_{S-J}, \a\>+M\<\rho^\vee_{S-J}, x y \i \a\><0.$$ By (b), $\a \in \Phi^+_J$. Therefore $v \in W_J$. Since $\l-M \rho^\vee_{S-J}$ and $\rho^\vee_{S-J}$ are dominant, $x y \i \a<0$. By the previous Lemma, we have that, $l(v \i y \i x)=l(y \i x)-l(v)$. Since $y \i x \in W^J$, $v \i y \i x \in W^J$. 

Now $\e^{-(\l-M \rho^\vee_{S-J}+y \i x M \rho^\vee_{S-J})} y \i x=v \e^{-\g} v \i y \i x$. We check that 

(c) $I(\g) \subsetneqq J$. 

Let $i \in S$ with $\<\g, \a_i\>=0$. Then $\<v \g, v \a_i\>=0$. Hence $v \a_i \in \Phi_J$ and $i \in J$. If $I(\g)=J$, then $\<v \g, \a\>=\<\g, v \i \a\>=0$ for any $\a \in \Phi_J^+$. By (b), $\<v \g, \a\> \ge 0$ for $\a \in \Phi^+$ and $v \g$ is dominant. Since $\g$ is dominant, $v \g=\g$. Notice that $v \in W^{I(\g)}$. Then $v=1$. That is a contradiction. (c) is proved. 
 
(d) $\<\g, \a_i\> \ge (\<\rho^\vee, \th\>+2)^{|I(\g)|+1}$ if $i \notin I(\g)$. 

If $i \notin J$, then $\<\g, \a_i\>=\<v \g, v \a_i\>$. We have that $v \a_i \in \Phi^+-\Phi^+_J$. By (b), $\<\g, \a_i\> \ge M$. If $i \in J-I(\g)$, then $\<\g, \a_i\>=\<v \g, v \a_i\>>0$. Notice that $\<v \g, v \a_i\>=\<\l-M \rho^\vee_{S-J}, v \a_i\>+M\<y \i x \rho^\vee_{S-J}, v \a_i\>$ and $\<\l-M \rho^\vee_{S-J}, v \a_i\>=0$. So $\<y \i x \rho^\vee_{S-J}, v \a_i\> \ge 1$ and $\<\g, \a_i\> \ge M$. (d) is proved. 

(e) $v \i y \i x v \in W^{I(\g)}$ and $\supp(v \i y \i x v)=S$. 

We have shown that $v \i y \i x \in W^J$ and $v \in W_J \cap W^{I(\g)}$. Then $v \i y \i x v \in W^{I(\g)}$ and $l(v \i y \i x v)=l(v \i y \i x)+l(v)$. So \begin{align*} \supp(v \i y \i x v) &=\supp(v \i y \i x) \cup \supp(v) \supset \supp(v \cdot v \i y \i x) \\ &=\supp(y \i x)=S.\end{align*} (e) is proved. 

By induction hypothesis for $(v, v \i y \i x, \g)$, $X_{v \e^{-\g} v \i y \i x}(1) \neq \emptyset$. By (a), $1 \in K \cdot_\s I \dot v \e^{-\g} \dot v \i \dot y \i \dot x I \subset K \cdot_\s I \dot x \e^{-\l} \dot y \i I$. Thus $X_{x \e^{-\l} y \i}(1) \neq \emptyset$ and the lemma holds in this case. \qed

\subsection{Proof of Theorem \ref{main2}}\label{main22} Assume that $[J, x, y] \cap Z_{J, w} \neq \emptyset$ for some $w \in W^J$ with $\supp(w)=S$. By Proposition \ref{key1}, $K \cdot_\s I \dot w \e^{-\l} I \subset K \cdot_\s I \dot x \e^{-\l} \dot y \i I$. By Proposition \ref{key2}, $1 \in K \cdot_\s I \dot w \e^{-\l} I$. Hence $1 \in K \cdot_\s I \dot x \e^{-\l} \dot y \i I$ and $X_{x \e^{-\l} y \i}(1) \neq \emptyset$.

\section{Further discussions}

In this section we discuss in more details the condition $[J, x, y] \nsubseteq \sqcup_{w \in W^J, \supp(w) \neq S} Z_{J, w}$.

\begin{lem}\label{e}
Let $J \subset S$, $x \in W^J$ and $y \in W$. If $\supp(y \i x) \neq S$, then $[J, x, y] \subset \sqcup_{w \in W^J, \supp(w) \neq S} Z_{J, w}$. 
\end{lem}

Suppose that $\supp(y \i x)=J' \neq S$. Let $\l \in Y^+$ with $I(\l)=J$. By Proposition \ref{key1}, it suffices to prove that $K \cdot I \dot x \e^{-\l} \dot y \i I \subset \sqcup_{w \in W^J, \supp(w) \neq S} K_{w \e^{-\l}}$. 

We may assume that $x=u_1 v_1$ and $y=u_2 v_2$ for $u_1, u_2 \in W^{J'}$ and $v_1, v_2 \in W_{J'}$. Since $x \in y W_{J'}$, we must have that $u_1=u_2$. Since $l(u_1 \i x)=l(v_1)=l(x)-l(u_1)$, we have that $v_1 \in W^J$. Then \begin{gather*} l(x \e^{-\l} y \i)=l(x \e^{-\l})+l(y)=l(x \e^{-\l})+l(v_2 \i)+l(u_2 \i)=l(x \e^{-\l} v_2 \i)+l(u_2 \i); \\ l(v_1 \e^{-\l} v_2 \i)=l(\e^{-\l} v_2 \i)-l(v_1)=l(\e^{-\l} v_2 \i)-l(x)+l(u_1)=l(u_1)+l(x \e^{-\l} v_2 \i). \end{gather*} Hence \begin{align*} K \cdot I \dot x \e^{-\l} \dot y \i I &=K \cdot I \dot x \e^{-\l} \dot v_2 \i I \dot u_2 \i I=K \cdot I \dot u_2 \i I \dot x \e^{-\l} \dot v_2 \i I \\ &=K \cdot I \dot v_1 \e^{-\l} \dot v_2 \i I.\end{align*} 

We prove by induction on $l(v_2)$ that \[\tag{a} K \cdot I \dot v_1 \e^{-\l} \dot v_2 \i I \subset \sqcup_{w \in W^J \cap W_{J'}} K_{w \e^{-\l}}.\]

Let $v_2=s_{i_1} \cdots s_{i_n}$ be a reduced expression. For $v_2=1$ this is clear. Assume now that $n>0$ and (a) holds for all $v'_2$ with $l(v'_2)<n$ but (a) fails for $v_2$. Then we have that \begin{align*} K \cdot I \dot v_1 \e^{-\l} \dot v_2 \i I &=K \cdot I \dot v_1 \e^{-\l} \dot v_2 \i \dot s_{i_1} I \dot s_{i_1} I=K \cdot I \dot s_{i_1} I \dot v_1 \e^{-\l} \dot v_2 \i \dot s_{i_1} I \\ & \subset K \cdot I \dot s_{i_1} \dot v_1 \e^{-\l} \dot v_2 \i \dot s_{i_1} I \cup K \cdot I \dot v_1 \e^{-\l} \dot v_2 \i \dot s_{i_1} I. \end{align*}
By induction hypothesis, $K \cdot I \dot v_1 \e^{-\l} \dot v_2 \i \dot s_{i_1} I\subset \sqcup_{w \in W^J \cap W_{J'}} K_{w \e^{-\l}}$. So $K \cdot I \dot s_{i_1} \dot v_1 \e^{-\l} \dot v_2 \i \dot s_{i_1} I \nsubseteq \sqcup_{w \in W^J \cap W_{J'}} K_{w \e^{-\l}}$. By induction hypothesis, $s_{i_1} v_1 \notin W^J$. Hence $s_{i_1} v_1=v_1 s_{j_1}$ for some $j_1 \in J$. Also $v_1, s_{i_1} \in W_{J'}$. So $j_1 \in J'$. Now $$K \cdot I \dot s_{i_1} \dot v_1 \e^{-\l} \dot v_2 \i \dot s_{i_1} I=K \cdot I \dot v_1 \e^{-\l} \dot s_{j_1} \dot v_2 \i \dot s_{i_1} I \nsubseteq \sqcup_{w \in W^J \cap W_{J'}} K_{w \e^{-\l}}.$$ Again by induction hypothesis, $s_{i_1} v_2 s_{j_1}>s_{i_1} v_2$ and $l(s_{i_1} v_2 s_{j_1})=n$. Now apply the same argument to $s_{i_1} v_2 s_{j_1}$ instead of $v_2$, we have that $s_{i_2} v_1=v_1 s_{j_2}$ for some $j_2 \in J \cap J'$ and $l(s_{i_2} s_{i_1} v_2 s_{j_1} s_{j_2})=n$. Repeat the same procedure, one can show that for $1 \le k \le n$ and $m \in \NN$, $s_{i_k} v_1^m=v_1^m s_j$ for some $j \in J \cap J'$. In other words, $\supp(v_2) \subset I(J, v_1)=I(S, v_1 \e^{-\l})$. Now by Lemma \ref{2}, $K \cdot I \dot v_1 \e^{-\l} \dot v_2 \i I \subset K_{v_1 \e^{-\l}}$. That is a contradiction. Therefore (a) always holds and the Lemma is proved. \qed

\begin{lem}\label{ne}
Let $J \subset S$, $x \in W^J$ and $y \in W$. If $\supp(u y \i x u \i)=S$ for all $u \in W_J$, then there exists $w \in W^J$ with $\supp(w)=S$ and $[J, x, y] \cap Z_{J, w} \neq \emptyset$. 
\end{lem}

By \cite[Corollary 2.6]{H4}, there exists $u \in W_J$, $w \in W^J$ and $v \in W_{I(J, w)}$ such that $u y \i x u \i=w v$. By our assumption, $\supp(w v)=S$. By the proof of \cite[Prop 3.2]{H5}, $\supp(w)=S$. Now \begin{align*} (\dot x, \dot y) \cdot h_J &=(\dot x \dot u \i, \dot y \dot u \i) \cdot h_J \in G_\D (\dot u \dot y \i \dot x \dot u \i, T) \cdot h_J \\ &=G_\D (\dot w \dot v, T) \cdot h_J \subset G_\D \cdot [J, w, v \i].\end{align*} By \cite[Proposition 1.10]{H3}, $(\dot x, \dot y) \cdot h_J \in Z_{J, w}$. The Lemma is proved. \qed

\begin{cor}\label{main3}
(1) Let $\l \in Y^+$ with $I(\l)=J \subsetneqq S$. Let $x \in W^J$ and $y \in W$. If $\supp(y \i x) \neq S$, then $X_{x \e^{-\l} y \i}(1)=\emptyset$.

(2) Let $J \subsetneqq S$, $x \in W^J$ and $y \in W$. Let $\l \in X^+$ with $I(\l)=J$ and $\<\l, \a_i\> \ge (\<\rho^\vee, \th\>+2)^{|J|+1}$. If $y=1$ and $\supp(x)=S$ or $\supp(u y \i x u \i)=S$ for all $u \in W_J$, then $X_{x \e^{-\l} y \i}(1) \neq \emptyset$. 
\end{cor}

\begin{rmk*}
Part (1) was first proved in \cite[Proposition 9.4.4]{GHKR2}. Here we give a different proof.  
\end{rmk*}

Part (1) follows from Proposition \ref{main1} and Lemma \ref{e}. The case where $y=1$ and $\supp(x)=S$ in part (2) follows from Proposition \ref{key2} and the case where $\supp(u y \i x u \i)=S$ for all $u \in W_J$ follows from Theorem \ref{main2} and Lemma \ref{ne}. \qed

\subsection{} In fact, the condition $\supp(u y \i x u \i)=S$ for all $u \in W_J$ can't be replaced by $\supp(y \i x)=S$ in general. The following is an example. 

Let $G=PGL_4$, $J=\{1\}$, $x=s_2 s_1 s_3 s_2$ and $y=s_3 s_2$. Then $\supp(y \i x)=S$ but $\supp(s_1 y \i x s_1)=\{2, 3\}$. In this case, one can show that $[\{1\}, x, y] \subset Z_{\{1\}, s_3 s_2}$. So $X_{x \e^{-\l} y \i}(1)=\emptyset$ for all $\l \in X^+$ with $I(\l)=J$. 

\section*{Acknowledgement} 
We thank George Lusztig and Ulrich G\"ortz for helpful discussions.


\end{document}